# Transformation of the linear difference equation into a system of the first order difference equations

M.I. Ayzatsky[1]

National Science Center Kharkov Institute of Physics and Technology (NSC KIPT),
610108, Kharkov, Ukraine

The transformation of the $N$th- order linear difference equation into a system of the first order difference equations is presented. The proposed transformation gives possibility to get new forms of the $N$-dimensional system of the first order equations that can be useful for analysis of the solutions of the $N$th- order difference equations. In particular, for the third-order linear difference equation the nonlinear second-order difference equation that plays the same role as the Riccati equation for second-order linear difference equation is obtained. The new form of the $N$-dimensional system of first order equations can be also used for finding the WKB solutions of the linear difference equation with coefficients that vary sufficiently slowly with index.

## 1 Introduction

It is common knowledge that a difference equation of order $N$

$$y_{k+N} + f_{N-1,k} y_{k+N-1} + f_{N-2,k} y_{k+N-2} + \ldots + f_{2,k} y_{k+2} + f_{1,k} y_{k+1} + f_{0,k} y_k + f_k = 0 \tag{1}$$

may be transformed in a standard way to a system of the $N$ first-order difference equations. For making such transformation we introduce a number of new variables (see, for example, [1,2])

$$x_k^{(i)} = y_{k+i-1}, \quad i = 1, 2, \ldots, N. \tag{2}$$

The difference equation (1) can be rewritten as

$$X_{k+1} = T_k X_k + F_k, \tag{3}$$

where $X_k = (x_{k+N-1}, x_{k+N-1}, \ldots, x_k)^T$, $F_k = (f_k, 0, \ldots, 0)^T$, and the companion matrix of (1) is

$$T_k = \begin{pmatrix} -f_{N-1,k} & -f_{N-2,k} & \ldots & -f_{1,k} & -f_{0,k} \\ 1 & 0 & \ldots & 0 & 0 \\ 0 & 1 & \ldots & 0 & 0 \\ & & & & \\ 0 & 0 & \ldots & 1 & 0 \end{pmatrix}. \tag{4}$$

In fact, we do not introduce new variables[2], we only do re-designations and still work with the elements of the same sequence $y_k$.

There is another kind of transformation that consists in representation of the solution $y_k$ of the equation (1) as the sum of the $N$ new unknown grid functions. For the case $N = 2$ it is discussed in [3,4]. By introducing $N$ new unknowns, instead of the one, we can impose $(N-1)$ additional conditions. Such approach gives new form of the $N$-dimensional system of first order equations, equivalent to the equation (1). In this note some generalization of the proposed transformation [3] is given. Analysis of literature [1,2,5,6,7,8,9,10,11,12,13,14,15,16,17,18,19,20] shows that it apparently has not been described earlier.

---

[1] M.I. Aizatskyi, N.I.Aizatsky; aizatsky@kipt.kharkov.ua

[2] In the case of transformation of differential equations we do introduce new variables $x^{(i)} = d^{i-1} y / dt^{i-1}, \quad i = 1, 2, \ldots, N$

## 2 Transformation the *N*th-order linear difference equation

We represent the solution of the difference equation (1) as the sum of new grid functions

$$y_k = \sum_{n=1}^{N} y_{n,k}. \tag{5}$$

By introducing $N$ new unknowns $y_{n,k}$ instead of the one $y_k$, we can impose additional conditions. These conditions we write in the form

$$y_{k+1} = \sum_{n=1}^{N} g_{1,n,k} y_{n,k}$$

$$y_{k+2} = \sum_{n=1}^{N} g_{2,n,k} y_{n,k} \tag{6}$$

$$....$$

$$y_{k+N-1} = \sum_{n=1}^{N} g_{N-1,n,k} y_{n,k}$$

where $g_{m,n,k}$ ($1 \leq n \leq N,\ 1 \leq m \leq N-1$) are the arbitrary sequences.

If

$$\det \begin{pmatrix} 1 & 1 & .. & 1 \\ g_{1,1,k} & g_{1,2,k} & \cdots & g_{1,N,k} \\ ... & ... & ... & ... \\ g_{N-1,1,k} & g_{N-1,2,k} & \cdots & g_{N-1,N,k} \end{pmatrix} \neq 0, \tag{7}$$

then the representation (5)-(6) is unique. Indeed, from (5) and (6) we can uniquely find $y_{n,k}$ as a linear combination of $y_k$. Using (1),(5) and (6) we can write such system of equations

$$y_{k+1} = \sum_{n=1}^{N} y_{n,k+1} = \sum_{n=1}^{N} g_{1,n,k} y_{n,k}$$

$$y_{k+2} = \sum_{n=1}^{N} g_{1,n,k+1} y_{n,k+1} = \sum_{n=1}^{N} g_{2,n,k} y_{n,k}$$

$$.... \tag{8}$$

$$y_{k+N-1} = \sum_{n=1}^{N} g_{N-2,n,k+1} y_{n,k+1} = \sum_{n=1}^{N} g_{N-1,n,k} y_{n,k}$$

$$\sum_{n=1}^{N} g_{N-1,n,k+1} y_{n,k+1} = -\sum_{n=1}^{N} \left( \sum_{m=1}^{N-1} f_{N-m,k} g_{N-m,n,k} + f_{0,k} \right) y_{n,k} - f_k$$

In matrix form

$$M_{k+1} Y_{k+1} = H_{k+1} Y_k + F_{k+1}, \tag{9}$$

where $Y_k = (y_{1,k}, y_{2,k}, ..., y_{N,k})^T$, $F_k = (0, 0..., -f_k)^T$,

$$M_{k+1} = \begin{pmatrix} 1 & 1 & .. & 1 \\ g_{1,1,k+1} & g_{1,2,k+1} & \cdots & g_{1,N,k+1} \\ ... & ... & ... & ... \\ g_{N-1,1,k+1} & g_{N-1,2,k+1} & \cdots & g_{N-1,N,k+1} \end{pmatrix}, \tag{10}$$

$$H_{k+1} = \begin{pmatrix} g_{1,1,k} & g_{1,2,k} & \cdots & g_{1,N,k} \\ \cdots & \cdots & \cdots & \cdots \\ g_{N-1,1,k} & g_{N-1,2,k} & \cdots & g_{N-1,N,k} \\ A_{1,k} & A_{2,k} & & A_{N,k} \end{pmatrix}, \quad (11)$$

$$A_{n,k} = -\left( \sum_{m=1}^{N-1} f_{N-m,k} g_{N-m,n,k} + f_{0,k} \right). \quad (12)$$

And finally we have the equation

$$Y_{k+1} = T_{k+1} Y_{k-1} + \overline{F}_k, \quad (13)$$

where

$$\begin{aligned} T_{k+1} &= M_{k+1}^{-1} H_{k+1} \\ \overline{F}_k &= M_{k+1}^{-1} F_{k+1} \end{aligned}. \quad (14)$$

As $g_{m,n,k}$ are the arbitrary sequences we can try to find such sequences that matrix $T_k$ will be a diagonal one. If it can be done, we find the solution of the difference equation (1).

If we chose $g_{m,n,k} = \left( \rho_k^{(n)} \right)^m$, where $\rho_k^{(n)}$ are the solutions of the characteristic equation

$$\rho_k^N + f_{N-1,k} \rho_k^{N-1} + f_{N-2,k} \rho_k^{N-2} + .. + f_{2,k} \rho_k^2 y_{k+2} + f_{1,k} \rho_k + f_{0,k} = 0, \quad (15)$$

the new form of the $N$-dimensional system of first order equations (13) can be used for finding the WKB solutions of the linear difference equation with coefficients that vary sufficiently slowly with index. For the case of the second-order difference equation it was proved in [4]. Bellow, we consider the case of the third-order difference equation and make generalization of the WKB solutions to the $N$th-order difference equation.

## 3 Transformation the third-order linear difference equation

Let's represent the solution of the linear third-order equation

$$y_{k+3} + f_{2,k} y_{k+2} + f_{1,k} y_{k+1} + f_{0,k} y_k + f_k = 0 \quad (16)$$

as the sum of three new functions

$$y_k = y_{1,k} + y_{2,k} + y_{3,k}. \quad (17)$$

Additional conditions we write in the form

$$\begin{aligned} y_{k+1} &= g_{1,1,k} y_{1,k} + g_{1,2,k} y_{2,k} + g_{1,3,k} y_{3,k} \\ y_{k+2} &= g_{2,1,k} y_{1,k} + g_{2,2,k} y_{2,k} + g_{2,3,k} y_{3,k} \end{aligned} \quad (18)$$

where $g_{m,n,k}$ ($m=1,2; n=1,2,3$) are the arbitrary sequences and

$$D_{k+1} = \begin{vmatrix} 1 & 1 & 1 \\ g_{1,1,k+1} & g_{1,2,k+1} & g_{1,3,k+1} \\ g_{2,1,k+1} & g_{2,2,k+1} & g_{2,3,k+1} \end{vmatrix} \neq 0. \quad (19)$$

Making the transformations that are given in section 2, we obtain a system of the first-order linear difference equations

$$y_{1,k+1}D_{k+1} = g_{1,1,k}y_{1,k}D_{k+1} +$$

$$+y_{1,k}\begin{bmatrix} g_{1,1,k}\{(g_{1,1,k+1}-g_{1,1,k})(g_{2,3,k+1}-g_{2,2,k+1})+(g_{2,1,k+1}-g_{2,1,k})(g_{1,2,k+1}-g_{1,3,k+1})\}+ \\ +x_{1,k}(g_{2,3,k+1}-g_{2,2,k+1})+(g_{1,2,k+1}-g_{1,3,k+1})x_{4,k} \end{bmatrix}$$

$$+y_{2,k}\begin{bmatrix} g_{1,2,k}\{(g_{1,2,k+1}-g_{1,2,k})(g_{2,3,k+1}-g_{2,2,k+1})+(g_{2,2,k+1}-g_{2,2,k})(g_{1,2,k+1}-g_{1,3,k+1})\}+ \\ +x_{2,k}(g_{2,3,k+1}-g_{2,2,k+1})+(g_{1,2,k+1}-g_{1,3,k+1})x_{5,k} \end{bmatrix} \quad (20)$$

$$+y_{3,k}\begin{bmatrix} g_{1,3,k}\{(g_{1,3,k+1}-g_{1,3,k})(g_{2,3,k+1}-g_{2,2,k+1})+(g_{2,3,k+1}-g_{2,3,k})(g_{1,2,k+1}-g_{1,3,k+1})\}+ \\ +x_{3,k}(g_{2,3,k+1}-g_{2,2,k+1})+(g_{1,2,k+1}-g_{1,3,k+1})x_{6,k} \end{bmatrix}$$

$$-f_k(g_{1,3,k+1}-g_{1,2,k+1})$$

$$y_{2,k+1}D_{k+1} =$$

$$+y_{1,k}\begin{bmatrix} g_{1,1,k}\{(g_{1,1,k+1}-g_{1,1,k})(g_{2,1,k+1}-g_{2,3,k+1})+(g_{2,1,k+1}-g_{2,1,k})(g_{1,3,k+1}-g_{1,1,k+1})\}+ \\ +x_{1,k}(g_{2,1,k+1}-g_{2,3,k+1})+(g_{1,3,k+1}-g_{1,1,k+1})x_{4,k} \end{bmatrix}$$

$$+y_{2,k}g_{1,2,k}D_{k+1} +$$

$$+y_{2,k}\begin{bmatrix} g_{1,2,k}\{(g_{1,2,k+1}-g_{1,2,k})(g_{2,1,k+1}-g_{2,3,k+1})+(g_{2,2,k+1}-g_{2,2,k})(g_{1,3,k+1}-g_{1,1,k+1})\}+ \\ +x_{2,k}(g_{2,1,k+1}-g_{2,3,k+1})+(g_{1,3,k+1}-g_{1,1,k+1})x_{5,k} \end{bmatrix} \quad (21)$$

$$+y_{3,k}\begin{bmatrix} g_{1,3,k}\{(g_{1,3,k+1}-g_{1,3,k})(g_{2,1,k+1}-g_{2,3,k+1})+(g_{2,3,k+1}-g_{2,3,k})(g_{1,3,k+1}-g_{1,1,k+1})\}+ \\ +x_{3,k}(g_{2,1,k+1}-g_{2,3,k+1})+(g_{1,3,k+1}-g_{1,1,k+1})x_{6,k} \end{bmatrix}$$

$$-f_k(g_{1,1,k+1}-g_{1,3,k+1})$$

$$y_{3,k+1}D_{k+1} =$$

$$y_{1,k}\begin{bmatrix} g_{1,1,k}\{(g_{1,1,k+1}-g_{1,1,k})(g_{2,2,k+1}-g_{2,1,k+1})+(g_{2,1,k+1}-g_{2,1,k})(g_{1,1,k+1}-g_{1,2,k+1})\}+ \\ +x_{1,k}(g_{2,2,k+1}-g_{2,1,k+1})+(g_{1,1,k+1}-g_{1,2,k+1})x_{4,k} \end{bmatrix}$$

$$+y_{2,k}\begin{bmatrix} g_{1,2,k}\{(g_{1,2,k+1}-g_{1,2,k})(g_{2,2,k+1}-g_{2,1,k+1})+(g_{2,2,k+1}-g_{2,2,k})(g_{1,1,k+1}-g_{1,2,k+1})+\} \\ +x_{2,k}(g_{2,2,k+1}-g_{2,1,k+1})+(g_{1,1,k+1}-g_{1,2,k+1})x_{5,k} \end{bmatrix} \quad (22)$$

$$+y_{3,k}g_{1,3,k}D_{k+1} +$$

$$+y_{3,k}\begin{bmatrix} g_{1,3,k}\{(g_{1,3,k+1}-g_{1,3,k})(g_{2,2,k+1}-g_{2,1,k+1})+(g_{2,3,k+1}-g_{2,3,k})(g_{1,1,k+1}-g_{1,2,k+1})\}+ \\ +x_{3,k}(g_{2,2,k+1}-g_{2,1,k+1})+(g_{1,1,k+1}-g_{1,2,k+1})x_{6,k} \end{bmatrix}$$

$$-f_k(g_{1,2,k+1}-g_{1,1,k+1})$$

where the following notations were introduced

$$x_{1,k} = (g_{1,1,k})^2 - g_{2,1,k}$$
$$x_{2,k} = (g_{1,2,k})^2 - g_{2,2,k}$$
$$x_{3,k} = (g_{1,3,k})^2 - g_{2,3,k}$$
$$x_{4,k} = g_{1,1,k}g_{2,1,k} + f_{2,k}g_{2,1,k} + f_{1,k}g_{1,1,k} + f_{0,k}$$
$$x_{5,k} = g_{1,2,k}g_{2,2,k} + f_{2,k}g_{2,2,k} + f_{1,k}g_{1,2,k} + f_{0,k}$$
$$x_{6,k} = g_{1,3,k}g_{2,3,k} + f_{2,k}g_{2,3,k} + f_{1,k}g_{1,3,k} + f_{0,k}$$

(23)

If we choose

$$g_{1,n,k} = \rho_k^{(n)},$$
$$g_{2,n,k} = \rho_k^{(n)2}, \; n = 1, 2, 3$$

(24)

where $\rho_k^{(n)}$ are the solutions of the equation

$$\rho_k^3 + f_{2,k}\rho_k^2 + f_{1,k}\rho_k + f_{0,k} = 0,$$

(25)

then $x_{i,k} = 0$, $i = 1,...,6$ and the system (20) - (22) takes the form

$$y_{1,k+1} = y_{1,k}\rho_k^{(1)} + y_{1,k}\rho_k^{(1)}\left(\rho_{k+1}^{(1)} - \rho_k^{(1)}\right)\left(\rho_{k+1}^{(3)} - \rho_{k+1}^{(2)}\right)\left[\left(\rho_{k+1}^{(3)} + \rho_{k+1}^{(2)}\right) - \left(\rho_{k+1}^{(1)} + \rho_k^{(1)}\right)\right]/D_{k+1} +$$
$$+ y_{2,k}\rho_k^{(2)}\left(\rho_{k+1}^{(2)} - \rho_k^{(2)}\right)\left(\rho_{k+1}^{(3)} - \rho_{k+1}^{(2)}\right)\left[\left(\rho_{k+1}^{(3)} + \rho_{k+1}^{(2)}\right) - \left(\rho_{k+1}^{(2)} + \rho_k^{(2)}\right)\right]/D_{k+1}$$
$$+ y_{3,k}\rho_k^{(3)}\left(\rho_{k+1}^{(3)} - \rho_k^{(3)}\right)\left(\rho_{k+1}^{(3)} - \rho_{k+1}^{(2)}\right)\left[\left(\rho_{k+1}^{(3)} + \rho_{k+1}^{(2)}\right) - \left(\rho_{k+1}^{(3)} + \rho_k^{(3)}\right)\right]/D_{k+1}$$
$$- f_k\left(\rho_{k+1}^{(3)} - \rho_{k+1}^{(2)}\right)/D_{k+1}$$

(26)

$$y_{2,k+1} = y_{1,k}\rho_k^{(1)}\left(\rho_{k+1}^{(1)} - \rho_k^{(1)}\right)\left(\rho_{k+1}^{(1)} - \rho_{k+1}^{(3)}\right)\left[\left(\rho_{k+1}^{(1)} + \rho_{k+1}^{(3)}\right) - \left(\rho_{k+1}^{(1)} + \rho_k^{(1)}\right)\right]/D_{k+1} +$$
$$+ y_{2,k}\rho_k^{(2)} + y_{2,k}\rho_k^{(2)}\left(\rho_{k+1}^{(2)} - \rho_k^{(2)}\right)\left(\rho_{k+1}^{(1)} - \rho_{k+1}^{(3)}\right)\left[\left(\rho_{k+1}^{(1)} + \rho_{k+1}^{(3)}\right) - \left(\rho_{k+1}^{(2)} + \rho_k^{(2)}\right)\right]/D_{k+1} +$$
$$+ y_{3,k}\rho_k^{(3)}\left(\rho_{k+1}^{(3)} - \rho_k^{(3)}\right)\left(\rho_{k+1}^{(1)} - \rho_{k+1}^{(3)}\right)\left[\left(\rho_{k+1}^{(1)} + \rho_{k+1}^{(3)}\right) - \left(\rho_{k+1}^{(3)} + \rho_k^{(3)}\right)\right]/D_{k+1} -$$
$$- f_k\left(\rho_{k+1}^{(1)} - \rho_{k+1}^{(3)}\right)/D_{k+1}$$

(27)

$$y_{3,k+1} = y_{1,k}\rho_k^{(1)}\left(\rho_{k+1}^{(1)} - \rho_k^{(1)}\right)\left(\rho_{k+1}^{(2)} - \rho_{k+1}^{(1)}\right)\left[\left(\rho_{k+1}^{(2)} + \rho_{k+1}^{(1)}\right) - \left(\rho_{k+1}^{(1)} + \rho_k^{(1)}\right)\right]/D_{k+1} +$$
$$+ y_{2,k}\rho_k^{(2)}\left(\rho_{k+1}^{(2)} - \rho_k^{(2)}\right)\left(\rho_{k+1}^{(2)} - \rho_{k+1}^{(1)}\right)\left[\left(\rho_{k+1}^{(2)} + \rho_{k+1}^{(1)}\right) - \left(\rho_{k+1}^{(2)} + \rho_k^{(2)}\right)\right]/D_{k+1} +$$
$$+ y_{3,k}\rho_k^{(3)} + y_{3,k}\rho_k^{(3)}\left(\rho_{k+1}^{(3)} - \rho_k^{(3)}\right)\left(\rho_{k+1}^{(2)} - \rho_{k+1}^{(1)}\right)\left[\left(\rho_{k+1}^{(2)} + \rho_{k+1}^{(1)}\right) - \left(\rho_{k+1}^{(3)} + \rho_k^{(3)}\right)\right]/D_{k+1} -$$
$$- f_k\left(\rho_{k+1}^{(2)} - \rho_{k+1}^{(1)}\right)/D_{k+1}$$

(28)

where $D_{k+1} = \left(\rho_{k+1}^{(2)} - \rho_{k+1}^{(1)}\right)\left(\rho_{k+1}^{(3)} - \rho_{k+1}^{(1)}\right)\left(\rho_{k+1}^{(3)} - \rho_{k+1}^{(2)}\right)$ - the Vandermonde determinant.

If the sequences $f_{0,k}, f_{1,k}, f_{2,k}$ vary sufficiently slowly with $k$ ($f_{0,k} = f_0(\varepsilon k)$, $f_{1,k} = f_1(\varepsilon k)$, $f_{2,k} = f_2(\varepsilon k)$, $0 \leq \varepsilon \ll 1$), then the differences $\left(\rho_{k+1}^{(n)} - \rho_k^{(n)}\right)$ are the small values and we can neglect the non-diagonal terms in the matrix $T_k$. This gives the WKB approximation

$$y_{1,k+1} \approx y_{1,k}\rho_k^{(1)} - y_{1,k}\rho_k^{(1)}\left(\rho_{k+1}^{(1)} - \rho_k^{(1)}\right)\left[\frac{1}{\left(\rho_{k+1}^{(1)} - \rho_{k+1}^{(2)}\right)} + \frac{1}{\left(\rho_{k+1}^{(1)} - \rho_{k+1}^{(3)}\right)}\right] - f_k\frac{\left(\rho_{k+1}^{(3)} - \rho_{k+1}^{(2)}\right)}{D_{k+1}},$$

(29)

$$y_{2,k+1} \approx y_{2,k}\rho_k^{(2)} - y_{2,k}\rho_k^{(2)}\left(\rho_{k+1}^{(2)} - \rho_k^{(2)}\right)\left[\frac{1}{\left(\rho_{k+1}^{(2)} - \rho_{k+1}^{(3)}\right)} + \frac{1}{\left(\rho_{k+1}^{(2)} - \rho_{k+1}^{(1)}\right)}\right] - f_k\frac{\left(\rho_{k+1}^{(1)} - \rho_{k+1}^{(3)}\right)}{D_{k+1}},$$

(30)

$$y_{3,k+1} \approx y_{3,k}\rho_k^{(3)} - y_{3,k}\rho_k^{(3)}\left(\rho_{k+1}^{(3)} - \rho_k^{(3)}\right)\left[\frac{1}{\left(\rho_{k+1}^{(3)} - \rho_{k+1}^{(1)}\right)} + \frac{1}{\left(\rho_{k+1}^{(3)} - \rho_{k+1}^{(2)}\right)}\right] - f_k \frac{\left(\rho_{k+1}^{(2)} - \rho_{k+1}^{(1)}\right)}{D_{k+1}}. \quad (31)$$

If we choose the sequences $g_{m,n,k}$ that are the solutions of the following equations

$$\begin{aligned}
x_{1,k} &= \left(g_{1,1,k}\right)^2 - g_{2,1,k} = -g_{1,1,k}\left(g_{1,1,k+1} - g_{1,1,k}\right), \\
x_{4,k} &= g_{1,1,k}g_{2,1,k} + f_{2,k}g_{2,1,k} + f_{1,k}g_{1,1,k} + f_{0,k} = -g_{1,1,k}\left(g_{2,1,k+1} - g_{2,1,k}\right), \\
x_{2,k} &= \left(g_{1,2,k}\right)^2 - g_{2,2,k} = -g_{1,2,k}\left(g_{1,2,k+1} - g_{1,2,k}\right), \\
x_{5,k} &= g_{1,2,k}g_{2,2,k} + f_{2,k}g_{2,2,k} + f_{1,k}g_{1,2,k} + f_{0,k} = -g_{1,2,k}\left(g_{2,2,k+1} - g_{2,2,k}\right), \\
x_{3,k} &= \left(g_{1,3,k}\right)^2 - g_{2,3,k} = -g_{1,3,k}\left(g_{1,3,k+1} - g_{1,3,k}\right), \\
x_{6,k} &= g_{1,3,k}g_{2,3,k} + f_{2,k}g_{2,3,k} + f_{1,k}g_{1,3,k} + f_{0,k} = -g_{1,3,k}\left(g_{2,3,k+1} - g_{2,3,k}\right),
\end{aligned} \quad (32)$$

the system (22) takes the form

$$\begin{aligned}
y_{1,k+1} &= y_{1,k}g_{1,1,k} - f_k \frac{\left(g_{1,3,k+1} - g_{1,2,k+1}\right)}{D_{k+1}}, \\
y_{2,k+1} &= y_{2,k}g_{1,2,k} - f_k \frac{\left(g_{1,1,k+1} - g_{1,3,k+1}\right)}{D_{k+1}}, \\
y_{3,k+1} &= y_{3,k}g_{1,3,k} - f_k \frac{\left(g_{1,2,k+1} - g_{1,1,k+1}\right)}{D_{k+1}}.
\end{aligned} \quad (33)$$

From (32) it follows that the sequences $g_{m,n,k}$ are the three different solutions of the system of the first-order nonlinear difference equations

$$\begin{aligned}
p_k^{(1)}\left(p_{k+1}^{(1)} - p_k^{(1)}\right) + p_k^{(1)2} - p_k^{(2)} &= 0 \\
p_k^{(1)}\left(p_{k+1}^{(2)} - p_k^{(2)}\right) + f_{2,k}p_k^{(2)} + f_{1,k}p_k^{(1)} + f_{0,k} + p_k^{(2)}p_k^{(1)} &= 0
\end{aligned}. \quad (34)$$

This system can be written as the second-order nonlinear difference equation

$$p_{k+2}^{(1)}p_{k+1}^{(1)}p_k^{(1)} + f_{2,k}p_{k+1}^{(1)}p_k^{(1)} + f_{1,k}p_k^{(1)} + f_{0,k} = 0. \quad (35)$$

For the third-order linear difference equation (16) the equation (35) (or system (34)) plays the same role as the Riccati equation for second-order linear difference equation.

The functions $y_{n,k} = \prod_k g_{1,n,k}$ are the linear independent ones, and the general solution of the homogeneous equation (16) ($f_k = 0$) is

$$y_k = \sum_{n=1}^{3} y_{n,k_0} \prod_{s=k_0}^{k-1} g_{1,n,s}. \quad (36)$$

There are other forms of system of the first order equations that can be obtained from the system (22) by choosing different sequences $g_{m,n,k}$.

Finding the WKB solutions of the linear difference equation (16) with coefficients that vary sufficiently slowly with index $k$ by finding the three iteration solutions of the equation (35) (compare with [21]) is not a simple procedure. So it seems preferable to use for the third-order difference equation the approach that leads us to the WKB equations (29)-(31).

## 4 The WKB approximation for the $N$th-order linear difference equation

Results of the work [4] and the section 3 show that the WKB equations for the $N$-order linear difference equation with coefficients that vary sufficiently slowly with index can be obtained by choosing of sequences $g_{m,n,k} = \left(\rho_k^{(n)}\right)^m$, where $\rho_k^{(n)}$ are the solutions of the characteristic equation (15) and taking into consideration only diagonal elements of the matrix $T_k$ in the equation

$$Y_{k+1} = T_{k+1} Y_k = M_{k+1}^{-1} H_{k+1} Y_k. \qquad (37)$$

If we chose $g_{m,n,k} = \left(\rho_k^{(n)}\right)^m$, the matrix $M_k$ transforms into the Vandermonde matrix and we can find its inverse matrix

$$M_{k,i,j}^{-1} = \frac{(-1)^{j-1} \sigma_{k,i}^{(N-j)}}{\prod_{\substack{s=1 \\ s \neq i}}^{N} \left(\rho_k^{(s)} - \rho_k^{(1)}\right)}, \qquad (38)$$

where $\sigma_{k,i}^{(j)} = \sum_{1 \leq m_1 < m_2 < \ldots < m_j \leq N} \prod_{s=1}^{j} \rho_k^{(m_s)} \left(1 - \delta_{m_s, i}\right)$. The matrix $H_{k+1}$ for such choice of sequences $g_{m,n,k}$ has the form

$$H_{k+1} = \begin{pmatrix} \rho_k^{(1)} & \rho_k^{(2)} & \ldots & \rho_k^{(N)} \\ \ldots & \ldots & \ldots & \ldots \\ \rho_k^{(1)N-1} & \rho_k^{(2)N-1} & \ldots & \rho_k^{(N)N-1} \\ \rho_k^{(1)N} & \rho_k^{(2)N} & \ldots & \rho_k^{(N)N} \end{pmatrix}. \qquad (39)$$

In the WKB approximation, we suppose that all elements of the matrix $M_{k+1}^{-1} H_{k+1}$ equals zero except the diagonal ones. In this case, the system of equations (37) can be rewritten as

$$y_{i,k+1} \approx y_{i,k} \sum_{j=1}^{N} M_{k+1,i,j}^{-1} \rho_k^{(i)j} = y_{i,k} \sum_{j=1}^{N} \left( \rho_k^{(i)j} \frac{(-1)^{j-1} \sigma_{k+1,i}^{(N-j)}}{\prod_{\substack{s=1 \\ s \neq i}}^{N} \left(\rho_{k+1}^{(s)} - \rho_{k+1}^{(1)}\right)} \right). \qquad (40)$$

These equations are generalization to the case of the of $N$ th-order difference equation the WKB equations obtained for the third (section 3) and second order difference equations [4].

## Conclusions

We presented transformations of the linear difference equation into a system of the first order difference equations. The proposed transformation gives possibility to get new forms of the $N$-dimensional system of first order equations that can be useful for analysis of the solutions of the $N$ th-order difference equation. In particular, for the third-order linear difference equation the nonlinear second-order difference equation that plays the same role as the Riccati equation for second-order linear equation is obtained. The new form of the $N$-dimensional system of first order equations can be also used for finding the WKB solutions of the linear difference equation with coefficients that vary sufficiently slowly with index.

## References

1 Elaydi S. An Introduction to Difference Equations. Springer, 2005

2 Cull P., Flahive M., Robson R. Difference Equations. From Rabbits to Chaos. Springer, 2005

3 M.I. Ayzatsky On the matrix form of second-order linear difference equations. http://lanl.arxiv.org/ftp/arxiv/papers//1703/1703.09608.pdf, LANL.arXiv.org e-print archives, 2017 4 M.I. Ayzatsky A note on the WKB solutions of difference equations. http://lanl.arxiv.org/ftp/arxiv/papers/1806/1806.02196.pdf, LANL.arXiv.org e-print archives, 2018

5 Jordan Charles. Calculus of Finite Differences. Chelsea Publishing Company, 1950

6 Levy H., Lessman F. Finite difference equations. Macmillan, 1961

7 Cull Paul, Flahive Mary, Robson Robby Difference Equations. From Rabbits to Chaos. Springer, 2005

8 Brand Louis. Differential and Difference Equations. John Wiley & Sons, 1966

9 Miller Kenneth S. An introduction to the calculus of finite differences and difference equations. Dover Publications, 1966

10 Wimp Jet. Computation with recurrence relations. Pitman Publishing, 1984

11 Samarskii A.A., Nikolaev E.S. Numerical Methods for Grid Equations. Volume I Direct Methods. Birkhäuser Basel (1989)

12 Kocic V.L., Ladas G. Global behavior of nonlinear difference equations of higher order with applications. Kluwer, 1993

13 Kulenovic Mustafa R.S., Ladas G. Dynamics of second order rational difference equations. Chapman and Hall_CRC, 2001

14 Lakshmikantham V., Trigiante Donato. Theory of Difference Equations. Numerical Methods and Applications. CRC Press, 2002

15 Sedaghat Hassan. Nonlinear Difference Equations. Theory with Applications to Social Science. Models-Springer Netherlands, 2003

16 Ashyralyev Allaberen , Sobolevskii Pavel E. New Difference Schemes for Partial Differential Equations.. Springer Basel AG, 2004

17 Grove E.A., Ladas G. Periodicities in nonlinear difference equations. CRC, 2005

18 Saber Elaydi. An Introduction to Difference Equations. Springer, 2005

19 Banasiak J. Mathematical Modelling in One Dimension. An Introduction via Difference and Differential Equations. Cambridge University Press, 2013

20 Mickens R.E. Difference Equations. Theory, Applications and Advanced Topics. Chapman & Hall_CRC, 2015

21 Braun P. A. WKB method for three-term recursion relations and quasienergies of an anharmonic oscillator, TMF, 1978,V.37, N.3, pp.355–370